%&amstex
\input amstex
\input amsppt.sty
\magnification=\magstep1
\hsize=30truecc
\baselineskip=16truept
\vsize=22.2truecm
\NoBlackBoxes
\nologo
\pageno=1
\topmatter
\TagsOnRight

\def\N{\Bbb N}
\def\Z{\Bbb Z}
\def\Q{\Bbb Q}

\def\l{\left}
\def\r{\right}
\def\b{\bigg}

\def\({\b(}
\def\[{\b[}
\def\){\b)}
\def\]{\b]}

\def\t{\text}
\def\f{\frac}
\def\mo{\roman{mod}}

\def\bi{\binom}
\def\eq{\equiv}

\def\ls{\leqslant}

\hbox{J. Number Theory 130(2010), no.\,4, 930--935.}
\medskip
\title Some curious congruences modulo primes\endtitle
\author Li-Lu Zhao and Zhi-Wei Sun \endauthor

\affil Department of Mathematics, Nanjing University
     \\Nanjing 210093, People's Republic of China
    \\ {\tt zhaolilu\@gmail.com},  \quad {\tt zwsun\@nju.edu.cn}
 \endaffil
\abstract Let $n$ be a positive odd integer and let $p>n+1$ be a prime.
We mainly derive the following congruence:
$$\sum_{0<i_1<\cdots<i_n<p}\l(\f{i_1}3\r)\f{(-1)^{i_1}}{i_1\cdots i_n}\eq0\ (\mo\ p).$$
\endabstract
\thanks 2010 {\it Mathematics Subject Classification}.
Primary 11A07; Secondary 11A41, 11T06.
\newline\indent
The second author is responsible for communications, and supported
by the National Natural Science Foundation (grant 10871087) of
China.
\endthanks
\endtopmatter
\document

\heading{1. Introduction}\endheading

Simple congruences modulo prime powers are of interest in number
theory. Here are some examples of such congruences:

 (a) (Wolstenholme) $\sum_{k=1}^{p-1}1/k\eq0\ (\mo\ p^2)$ for any prime $p>3$.
 \medskip

 (b) (Z. W. Sun [S02, (1.13)]) For each
 prime $p>3$ we have
  $$\sum_{0<k<p/2}\f{3^k}k\eq\sum_{0<k<p/6}\f{(-1)^k}k\ (\mo\ p).$$

 (c) (Z. W. Sun [S07, Theorem 1.2]) If $p$ is a prime and $a,n\in\N=\{0,1,2,\ldots\}$, then
 $$\f1{\lfloor n/p^a\rfloor!}\sum_{k\eq0\,(\mo\ p^a)}(-1)^k\bi nk\l(-\f k{p^a}\r)^{\lfloor n/p^a\rfloor}\eq1\ (\mo\ p).$$

 (d) (Z. W. Sun and R. Tauraso [ST, Corollary 1.1]) For any prime $p$ and $a\in\Z^+=\{1,2,3,\ldots\}$ we have
 $$\sum_{k=0}^{p^a-1}\bi{2k}k\eq\l(\f{p^a}3\r)\ (\mo\ p^2),$$
where $(\f{\cdot}3)$ is the Legendre symbol.

\medskip
Let $p>3$ be a prime. In 2008, during his study of
$\sum_{k=0}^{p-1}\binom{2k}k$ modulo powers of $p$ with R. Tauraso,
the second author conjectured that
$$\sum_{0<i<j<k<p}\(\f{i}{3}\)\f{(-1)^i}{ijk}\eq 0\ (\mo\ p),\tag1.1$$
i.e.,
$$\sum\Sb0<i<j<k<p\\i\eq 1,2\,(\mo\ 6)\endSb\f{1}{ijk}\eq \sum\Sb 0<i<j<k<p\\i\eq 4,5\,(\mo\ 6)\endSb
\f{1}{ijk}\ (\mo\ p).\tag 1.2$$

In this paper we confirm the above conjecture of Sun by establishing the
following general theorem.

 \proclaim{Theorem 1.1} Let $n\in \Z^+$ and let $p>n+1$ be a prime.

{\rm (i)} If $n$ is odd, then
$$\sum\Sb0<i_1<\cdots<i_n<p\\i_1\eq 1,2\,(\mo\ 6)\endSb\f{1}{i_1\cdots i_n}
\eq \sum\Sb 0<i_1<\cdots<i_n<p\\i_1\eq 4,5\,(\mo\
6)\endSb\f{1}{i_1\cdots i_n}\ (\mo\ p).\tag 1.3$$

{\rm (i)} If $n$ is even, then
$$\sum\Sb0<i_1<\cdots<i_n<p\\i_1\eq 0\,(\mo\ 3)\endSb\f{(-1)^{i_1}}{i_1\cdots i_n}
\eq 2\sum\Sb 0<i_1<\cdots<i_n<p\\i_1\eq 2,3,4\,(\mo\
6)\endSb\f{1}{i_1\cdots i_n}\ (\mo\ p).\tag 1.4$$
\endproclaim

We deduce Theorem 1.1 from our following result.

\proclaim{Theorem 1.2} Let $n\in \Z^+$ and let $p>n+1$ be a prime.
 Set
 $$F_n(x)=\sum_{0<i_1<\dots <i_n<p}\f{x^{i_1}}{i_1\cdots i_n}\in\Z_p[x],\tag1.5$$
 where $\Z_p$ denotes the integral ring of the $p$-adic field $\Q_p$.
Then we have
$$F_n(1-x)\eq(-1)^{n-1}F_n(x)\ (\mo\ p),\tag1.6$$
i.e., all the coefficients of $F_n(1-x)-(-1)^{n-1}F_n(x)$ are congruent to $0$ modulo $p$.
\endproclaim

In the next section we use Theorem 1.2 to prove Theorem 1.1.
Section 3 is devoted to our proof of Theorem 1.2.

\heading{2. Theorem 1.2 implies Theorem 1.1}\endheading

\medskip
\noindent{\it Proof of Theorem 1.1 via Theorem 1.2}. (1.3) holds trivially when $p=3$ and $n=1$.
Below we assume that $p>3$.

Let $\omega$ be a primitive cubic root of unity in an extension field over $\Q_p$.
Then, in the ring $\Z_p[\omega]$ we have the congruence
$$F_n(-\omega^2)=F_n(1+\omega)\eq(-1)^{n-1}F_n(-\omega)\ (\mo\ p).\tag2.1$$

For $r\in\Z$ we set
$$S_r=\sum\Sb0<i_1<\cdots<i_n<p\\i_1\eq r\,(\mo\ 6)\endSb\f1{i_1\cdots i_n}.$$
Clearly
$$\align F_n(-\omega)=&S_0-\omega S_1+\omega^2 S_2-S_3+\omega S_4-\omega^2 S_5
\\=&S_0-S_3-\omega(S_1-S_4)+(-1-\omega)(S_2-S_5)
\\=&S_0-S_3-S_2+S_5-\omega(S_1+S_2-S_4-S_5).
\endalign$$
Similarly,
$$ F_n(-\omega^2)=S_0-S_3-S_2+S_5-\omega^2(S_1+S_2-S_4-S_5).$$
Thus
$$\align F_n(-\omega)+F_n(-\omega^2)=&2(S_0-S_2-S_3+S_5)+S_1+S_2-S_4-S_5
\\=&2S_0+S_1-S_2-2S_3-S_4+S_5
\endalign$$
and
$$F_n(-\omega)-F_n(-\omega^2)=(\omega^2-\omega)(S_1+S_2-S_4-S_5).$$
Note that $(\omega-1)(\omega^2-1)=3$ is relatively prime to $p$. Therefore, by (2.1),
if $2\nmid n$ then $$S_1+S_2-S_4-S_5\eq0\ (\mo\ p);\tag2.2$$
if $2\mid n$ then
$$2S_0+S_1-S_2-2S_3-S_4+S_5\eq0\ (\mo\ p).\tag2.3$$
To conclude the proof we only need to show that (2.3) is equivalent to
$$S_0-S_3\eq 2(S_2+S_3+S_4)\ (\mo\ p).\tag2.4$$

Recall that
$$x^{p-1}-1\eq\prod_{j=1}^{p-1}(x-j)\eq\prod_{i=1}^{p-1}\l(x-\f1i\r)\ (\mo\ p)$$
(cf. Proposition 4.1.1 of [IR, p.\,40]). Comparing the coefficients
of $x^{p-1-n}$ we get that
$$\sum_{0<i_1<\cdots<i_n<p}\f{1}{i_1\cdots i_n}\eq 0\ (\mo\ p).\tag 2.5$$
So $\sum_{r=0}^5 S_r\eq0\ (\mo\ p)$, which implies the equivalence of (2.3) and (2.4). We are done. \qed

\heading{3. Proof of Theorem 1.2}\endheading

\medskip
\noindent{\it Proof of Theorem 1.2}. We use induction on $n$.

Observe that
$$\sum_{i=1}^{p-1}\bi pi(-1)^{i-1}x^i=1+(-x)^p-\sum_{i=0}^p\bi pi(-x)^i=1-x^p-(1-x)^p.$$
For $i=1,\ldots,p-1$ clearly
$$\f{(-1)^{i-1}}p\bi pi=\f{(-1)^{i-1}}i\bi{p-1}{i-1}\eq\f1i\ (\mo\ p).$$
Thus
$$F_1(x)\eq\f1p\sum_{i=1}^{p-1}\bi pi(-1)^{i-1}x^i=\f{1-x^p-(1-x)^p}p\ (\mo\ p)$$
and hence $F_1(1-x)\eq F_1(x)\ (\mo\ p)$ as desired. This proves
(1.6) for $n=1$.

For the induction step we need to do some preparation. For
$$P(x)=\sum_{i=0}^ma_ix^i\in\Z_p[x],$$
we define its {\it formal derivative} by
$$\f {\t d}{{\t d}x}P(x)=\sum_{0<i\ls m}ia_ix^{i-1}.$$
If $1\ls m\ls p-1$ and $\f{\t d}{{\t d}x} P(x)\eq0\ (\mo\ p)$, then
$a_i\eq0\ (\mo\ p)$ for all $i=1,\ldots,m$, and hence $P(x)\eq
a_0=P(0)\ (\mo\ p)$.

Now assume that $1<n<p-1$ and $F_{n-1}(1-x)\eq (-1)^{n-2}F_{n-1}(x)\
(\mo\ p)$. Then
$$\align \f {\t d}{{\t d}x}F_n(x)=&\sum_{0<i_1<\dots <i_n<p}\f{x^{i_1-1}}{i_2\cdots i_n}
=\sum_{1<i_2<\dots <i_n<p}\f{1}{i_2\cdots i_n}\sum_{i_1=1}^{i_2-1}x^{i_1-1}
\\=&\sum_{0<i_2<\dots <i_n<p}\f{1}{i_2\cdots i_n}\cdot
\f{x^{i_2-1}-1}{x-1}
\\=&\f{F_{n-1}(x)}{x(x-1)}-\f1{x-1}\sum_{0<i_2<\dots <i_n<p}\f{1}{i_2\cdots i_n}
\endalign$$
and hence
$$\align &\f {\t d}{{\t d}x}\l(F_n(1-x)-(-1)^{n-1}F_n(x)\r)
\\=&-\(\f{F_{n-1}(1-x)}{(1-x)(1-x-1)}-\f1{(1-x)-1}\sum_{0<i_2<\cdots<i_n<p}\f1{i_2\cdots i_n}\)
\\&+(-1)^n\(\f{F_{n-1}(x)}{x(x-1)}-\f1{x-1}\sum_{0<i_2<\cdots<i_n<p}\f1{i_2\cdots i_n}\)
\\=&\f{(-1)^nF_{n-1}(x)-F_{n-1}(1-x)}{x(x-1)}-\l(\f1x+\f{(-1)^n}{x-1}\r)\sum_{0<i_2<\cdots<i_n<p}\f1{i_2\cdots i_n}.
\endalign$$
Combining this with the induction hypothesis and (2.5), we obtain
$$x(x-1)\f {\t d}{{\t d}x}(F_n(1-x)-(-1)^{n-1}F_n(x))\eq0\ (\mo\ p).$$
For the finite field $\Bbb F_p=\Z/p\Z$, it is well known that $\Bbb
F_p[x]$ is a principal ideal domain. So we have
$$\f {\t d}{{\t d}x}(F_n(1-x)-(-1)^{n-1}F_n(x))\eq0\ (\mo\ p)$$
and hence
$$\align &F_n(1-x)-(-1)^{n-1}F_n(x)
\\\eq &F_n(1)+(-1)^nF_n(0)=\sum_{0<i_1<\cdots<i_n<p}\f1{i_1\cdots i_n}\eq0\ (\mo\ p)
\endalign$$
with the help of (2.5). This concludes the induction step and we are done. \qed


\widestnumber\key{S07}

\Refs

\ref\key IR\by K. Ireland and M. Rosen \book A Classical
Introduction to Modern Number Theory \publ 2nd ed., Grad. Texts in
Math. 84, Springer, New York, 1990\endref

\ref\key S02\by Z. W. Sun\paper On the sum $\sum_{k\eq r\,(\mo\
m)}\bi nk$ and related congruences\jour Israel J. Math.\vol 128\yr
2002\pages 135--156\endref

\ref\key S07\by Z. W. Sun\paper Combinatorial congruences and
Stirling numbers\jour Acta Arith.\vol 126\yr 2007\pages
387--398\endref

\ref\key ST\by Z. W. Sun and R. Tauraso\paper On some new
congruences for binomial coefficients \jour preprint,
arXiv:0709.1665\endref

\endRefs
\enddocument